\input amstex
\documentstyle{amsppt}

\tolerance=500

\def \cC{\Cal C}
\def \cF{\Cal F}
\def \cI{\Cal I}
\def \cO{\Cal O}
\def \cK{\Cal K}
\def \cN{\Cal N}

\def \bP{\bold P}
\def \bH{\bold H}
\def \bZ{\bold Z}
\def \bN{\bold N}

\topmatter

\title Arithmetically Cohen--Macaulay curves in $\bP^4$\\
of degree $4$ and genus $0$
\endtitle
\rightheadtext{Arithmetically Cohen--Macaucaulay curves}
\author Mireille Martin-Deschamps and Ragni Piene \endauthor
\leftheadtext{M. Martin-Deschamps and R. Piene}
\address
D\'epartement de Math\'ematiques et d'Informatique, \'Ecole Normale
Sup\'e-rieure, 45
rue d'Ulm, F-75230 Paris C\'edex 05, France
\endaddress
\email Mireille.Deschamps\@ens.fr \endemail
\address Matematisk institutt, Universitetet i Oslo, P.B. 1053 Blindern,
N-0316 Oslo, Norway \endaddress
\email ragnip\@math.uio.no \endemail
\thanks The first author was supported in part by the Norwegian
Research Council (Matematisk Seminar). The second author was
supported in part by \'Ecole Normale Sup\'erieure and the Nansen Fund.
\endthanks
\abstract
We show that the arithmetically Cohen--Macaulay (ACM) curves of degree
$4$ and genus $0$ in $\bP^4$ form an irreducible subset of the Hilbert scheme.
Using this, we show that the singular locus of the corresponding component
of the Hilbert scheme has dimension greater than $6$.
Moreover, we describe the structures of all
ACM curves of Hilb$^{4m+1}(\bP^4)$.
\endabstract
\endtopmatter
\document
\head{Introduction}\endhead
Let $k$ be an algebraically closed field, and let $\bP^n=\bP^n_k$ denote
projective $n$-space over $k$. If $C\subset
\bP^n$ is a subscheme, we denote by
$\cI_C$ its sheaf of ideals, and by $I_C$ its homogeneous saturated
 ideal
in the ring $S=k[X_0, \ldots, X_n]$. We recall that $C$ is
{\it locally Cohen-Macaulay} if all its local rings are Cohen--Macaulay
(i.e., if their depth is equal to their dimension) and that $C$
is
{\it arithmetically
Cohen--Macaulay} (ACM) if its homogeneous coordinate ring
$S/I_C$ is Cohen--Macaulay. If $C$ is ACM, then
it is locally Cohen--Macaulay. By a {\it curve} we
shall mean a
locally Cohen--Macaulay scheme of pure dimension $1$. Note that an
 ACM
curve $C\subset \bP^n$ is connected.

If $P(m)$ is a polynomial, we denote
by $Hilb^{P(m)}(\bP^n)$ the Hilbert scheme parametrizing subschemes
 of
$\bP^n$ with Hilbert polynomial $P(m)$. A $1$-dimensional subscheme
$C\subset \bP^n$ of degree $d$ and (arithmetic) genus $g$ has Hilbert
polynomial $dm+1-g$.
We shall say that a curve of degree
$d$ and arithmetic genus $g$ is a {\it curve of type} $(d,g)$ (or, for short,
a $(d,g)$).
The points of $Hilb^{P(m)}(\bP^n)$
corresponding to ACM curves form an open subset.

A {\it rational normal (RN) curve} in $\bP^n$ is a
$\bP^1$ embedded in $\bP^n$ by the sections of $\cO_{\bP^1}(n)$.
The normal bundle of a RN curve is $\cO_{\bP^1}(n+2)^{\oplus
{n-1}}$ (see e.g. [S]), and its
Hilbert polynomial is $P(m)=nm+1$. The RN curves correspond to the points
of a smooth,
irreducible open subscheme of $Hilb^{P(m)}(\bP^n)$ of dimension
$(n-1)(n+3)$. Moreover, the RN curves are the only reduced and
irreducible ACM curves with Hilbert polynomial $nm+1$.

For $n=3$, the scheme $Hilb^{3m+1}(\bP^3)$ was studied in [P--S], where it was
shown
that $Hilb^{3m+1}(\bP^3)$ consists of two smooth components,
intersecting transversally. The component containing the points
corresponding to RN curves of degree 3 (the {\it twisted
 cubics})
has dimension 12, and the other component --- consisting of plane cubics with an
isolated or embedded point --- has dimension 15. Moreover, in this
 case,
every locally Cohen-Macaulay curve is also ACM ([E], [H2]
Prop.3.5, [P--S]),
and the saturated ideal of such a curve is generated by three
independent quadratic forms. The fact that the component containing
the twisted cubics is smooth, has been useful in enumerative
questions concerning twisted cubics.

Except for the first fundamental results concerning existence and
basic
properties of Hilbert schemes ([G], [H1]), there are very
few results describing Hilbert schemes in general. A. Reeves has
shown that the diameter of the Hilbert scheme $Hilb^{P(m)}(\bP^n)$ is $\leq
2d+2$,
where $d$ is the degree of $P(n)$ ([R], Thm. 8) --- hence, in the case of curves
($d=1$) the diameter is always $\leq 4$. She also
describes points on
various components, but unfortunately this is not sufficient
to describe the different components, nor even to determine their number.

For $n\geq 4$ the scheme
$Hilb^{nm+1}(\bP^n)$ has more than two components. For example, in addition
to the
component containing the RN curves, there is a component
--- also of dimension $(n-1)(n+3)$ ---
with general
point corresponding to an elliptic curve of degree $n-1$ (in a $\bP^{n-2}$)
union a disjoint line;
moreover, there is a component --- of dimension ${1\over 2}
(n^3-2n^2+11n-12)$ ---
with general point
corresponding to a plane smooth curve of degree $n$ union ${1\over 2}(n-1)(n-2)$
points,
and so on. There seems to be no reason to expect the components
to be nonsingular.

The easiest way to exhibit singularities of a component, is to find a curve
such that the corresponding point on the Hilbert scheme
belongs to only one component and such that
the tangent space to the Hilbert scheme at the point is greater than the
dimension of the component. When the ideal of the curve is given, it is often
possible to compute the tangent space, but it is not always easy to determine
the number of components containing the curve.

For $n\geq 4$ we shall see that the component of $Hilb^{nm+1}(\bP^n)$
containing the RN curves contains a large subvariety, isomorphic to
the Grassmann variety Grass(1,n) of lines in $\bP^n$, where the tangent
space is ``too big''. These
points correspond to an ACM $n$-fold structure on a line. Since the set of ACM
curves is open, in order to show that the component is singular at these points,
it suffices to show that the set of ACM curves is irreducible.

It is reasonable to believe that the set of ACM curves is
irreducible if (and only if) $n\leq 7$, and hence that,
in this case, any ACM curve can be smoothed (to a RN
curve). For $n\geq 8$, however, the situation is different:
the lowest length example of a fat point which is not smoothable is a
point of degree
$8$ in $4$-dimensional space, such that the ideal of the point is generated by
seven quadratic forms (see [I], p.310). This point can be reembedded in $\bP^7$
such that the cone over the point is an
ACM curve (of genus 0) in $\bP^8$ which is not smoothable ([C2]).

The aim of this paper is to show that in the case
$n=4$, the set of ACM curves is indeed irreducible. In particular,
the possible {\it non-reduced} structures of an ACM curve with Hilbert
polynomial $4m+1$ are limited, and we shall describe the possible structures.

In the first section we gather some known facts about ACM curves in $\bP^n$
of degree $n$ and genus $0$.
In the next section we restrict to the case $n=4$. The strategy for proving that
the set of ACM curves is irreducible, is to show that any ACM curve is the limit
of a family of (possibly degenerated) RN curves, and to do this by
reducing to the case of curves of smaller degree (in a space of
smaller dimension). The difficult cases are the
non-reduced curves. When the underlying reduced curve contains a line, we
``remove'' the line by intersecting with a hyperplane, to get a
degeneration of a RN curve of degree $3$.
In certain cases, we need to use projections onto a hyperplane ---
the method is explained in Section 3.

As a byproduct of the methods used to prove irreducibility of the ACM
curves, we
obtain a fairly explicit description of the possible ACM curves in
$Hilb^{4m+1}(\bP^4)$. These results are gathered in Section 4.
\medskip
\noindent{\bf Acknowlegdment.} We would like to thank Jan A. Christophersen
for fruitful discussions.

\head{1. ACM curves in $\bP^n$ with Hilbert polynomial $nm+1$}\endhead
The ACM curves in $\bP^n$ with Hilbert polynomial $nm+1$ have been studied
by various people ([W], [E--R--S], [C1]). We shall recall some
results.

\proclaim{Proposition 1.1} Let
$C\subset \bP^n$ be an ACM
curve with Hilbert polynomial $nm+1$. Then its homogeneous ideal
$I_C$ is generated by ${1\over 2}n(n-1)$ independent quadratic forms.
\endproclaim

\demo{Proof} A connected curve $C$ of genus $0$ is ACM if and
only if $H^1(\bP^n,\cI_C(m))=0$ for all $m\geq0$.
Therefore, the sheaf of ideals
$\cI_C$ is $2$-regular. By the Castelnuovo--Mumford
regularity criterion it follows that $I_C$ is generated by quadrics,
and the number of independent quadrics is $$\text{dim}~H^0(C,\cI_C(2))=
{{n+2}\choose 2}-(2n+1)={1\over 2}n(n-1).~~\qed$$\enddemo

\medbreak
\subhead{1.2 The n-fold line}\endsubhead
By an {\it n-fold line} we shall mean any curve projectively
equivalent to $L_n$, the line
$X_1=\cdots =X_{n-1}=0$ with $n$-fold stucture defined by
$I_{L_n}=(X_1, \ldots, X_{n-1})^2$. The curve $L_n$ is ACM, and
there is a
$1$-parameter deformation of $L_n$ to a RN curve: consider the matrix
$$\pmatrix tX_0&X_1&\ldots &X_{n-2}&X_{n-1}\\
X_1&X_2&\ldots &X_{n-1}&tX_n
\endpmatrix
$$
For $t\neq 0$ the vanishing of the $(2\times 2)$-minors defines a
RN curve,
whereas for $t=0$, the minors generate the ideal of $L_n$.

The $n$-fold line is a particular multiple structure of degree n on a
line. It can be characterized in the following way:

\proclaim{Proposition 1.3} Let $C\subset \bP^n$ be an ACM curve of degree $n$
and genus $0$ which  contains  a line $L$. If the Zariski tangent space
$T_{C,x}$ to $C$ at $x$ has dimension $n$
for all $x\in L$, then $C$ is an n-fold line.\endproclaim

\demo{Proof} The dimension of the Zariski tangent space
$T_{C,x}$
is equal to the dimension of the projective tangent space
to $C$ at $x$, and
the latter is equal to the dimension of the null space of the Jacobian matrix
of $C$ at $x$.
We may assume that $I=(X_1,\cdots ,X_{n-1})$ is the ideal of the line. We know
that $I_C$ is generated by ${1\over 2}n(n-1)$ independent quadratic
forms belonging to $I$. All their partial derivatives with respect
to $X_1,\cdots ,X_{n-1}$ are also in
$I$, so $I_C$ is contained in $I^2$, therefore $I_C=I^2$ for degree reasons.
\qed \enddemo

\proclaim{Proposition 1.4} If
$C\subset \bP^n$ is an ACM
curve with Hilbert polynomial $nm+1$, then $C$ specializes to an $n$-fold line.
\endproclaim

\demo{Proof} ([C1], [E--R--S] Ex.4.1, [W] Rem.2.2.1) The quadratic generators
of the homogeneous ideal $I_C$ can be chosen in a
particular way: $C$ is ACM if and only if
$H^1(\bP^n,\cI_C(m))=0$ for all $m\geq0$, hence if and
only if dim$(S/I_C)_m=nm+1$ for all $m\geq0$. Choose a regular sequence
$L,M\in S_1$ for $S/I_C$ and put $A=S/(I_C+(L,M))$. One computes
dim$A_0=1$ and dim$A_1=n-1$, while dim$A_m=0$ for
$m\geq2$. By a coordinate change, we may assume $L=X_0$ and $M=X_n$, so
that $A$ is a quotient of $k[X_1, \ldots, X_{n-1}]$. But
then $A$ must be $k[X_1, \ldots, X_{n-1}]/(X_1, \ldots, X_{n-1})^2$.
It follows that the generators of $I_C$ may be written $$X_iX_j-X_0L_{i,j}
-X_nM_{i,j}-Q_{i,j},$$ where $1\leq i\leq j\leq n-1$, the $L_{i,j}$ and
$M_{i,j}$ are linear forms in $k[X_1, \ldots, X_{n-1}]$ and the $Q_{i,j}$
are quadratic forms in
$k[X_0,X_n]$, subject to certain conditions ensuring that $C$ is
an ACM curve. Finally one checks that setting $L_{i,j}=M_{i,j}=Q_{i,j}=0$
gives a flat specialization. \qed \enddemo

\medskip
Let $L\subset \bP^n$ be a n-fold line, with normal bundle ${\cN}_L$.
Christophersen ([C1]) showed that $H^0(L,{\cN}_L)=
\text{Hom}(I_L/I_L^2,S/I_L)_0$ and has dimension $n(n-1)^2$. Since $L$ lies
on the component of $Hilb^{nm+1}(\bP^n)$ containing the RN curves, which has
dimension $(n-1)(n+3)$, it follows that $L$ is a singular point on
$Hilb^{nm+1}(\bP^n)$ whenever $n(n-1)^2 \geq (n-1)(n+3)$, i.e., when $n\geq 4$.
In order to determine whether the $n$-fold lines also are singular points on
the {\it component}, one needs to check that they do not lie on any {\it other}
component of the Hilbert scheme. We shall prove (in Section 3)
that this holds for $n=4$.

\proclaim{Proposition 1.5} {\rm ([W], Prop. 5.6, p. 257)} Let $C\subset \bP^n$
 be a reduced curve with
Hilbert polynomial $nm+1$, not contained in a hyperplane. Then $C$
can be smoothed in $\bP^n$.\endproclaim

Note that it follows from this proposition, since any ACM curve is connected,
that a reduced ACM curve with Hilbert polynomial $nm+1$ can be smoothed
to a RN curve.

\head{2. The case $n = 4$}\endhead
Let $\bH^0$ denote the open subscheme of
$\bH=Hilb^{4m+1}(\bP^4)$ consisting of ACM curves.
We shall prove that $\bH^0$ is contained in the  closure $\bH'$
of the open subscheme consisting of the RN curves in $\bH$.

Applying Propositions 1.1 and 1.5 to the case $n=4$, we get that for
 any $C\in
\bH^0$, the ideal $I_C$ is generated by 6 independent quadratic
forms, and if $C\in \bH^0$ is reduced, then $C\in\bH'$.

When $C$ is ACM and not reduced, we
have  the following possibilities ($C_{red}$ is the underlying
reduced curve of $C$):

Case I: deg $C_{red}=3$. Then $C$ is the union of a double line
 and
a reduced curve of degree 2.

Case II: deg $C_{red}=2$. Then $C_{red}$ is a plane conic,
and we have 2 subcases:

-- II.1: $C$ is a double structure on a reduced conic;

-- II.2: $C$ is  the union of a reduced  line and a triple structure
on another line.

Case III: deg $C_{red}=1$. Then $C$ is a quadruple structure
on a line.

In each case, there exists a hyperplane such that its intersection
with $C$ contains a curve. This will give us, at least when the residual
intersection is a line, a way of constructing the curve as a limit
 of a
family of reduced ACM curves.

We shall need a technical result.

\proclaim{Lemma 2.1} Let $L$, $L_1$, $L_2$, $L_3$ be four independent linear
 forms,
and $Q_1$, $Q_2$, $Q_3$ three quadratic forms contained in the ideal
$(L_1,L_2,L_3)$. Then we have
$$(L,Q_1,Q_2,Q_3)\cap (L_1,L_2,L_3)=(LL_1,LL_2,LL_3,Q_1,Q_2,Q_3).$$
\endproclaim
\demo{Proof} $\alpha L+\sum \beta_iQ_i\in (L_1,L_2,L_3)$ if and only if
 $\alpha
\in (L_1,L_2,L_3)$.\qed \enddemo

\medskip

\proclaim{Proposition 2.2} Let $C\in \bH^0$ and let $L\in S$ be a linear
form such that the
corresponding  hyperplane $H_L$ intersects $C$ in a curve. Denote by $C'$
the largest curve contained in $C\cap H_L$. We have an exact sequence
 of
graded algebras:
$$0@>>> S/I_{\Gamma}(-1) @>{\cdot L}>> S/I_C @>>> S/I_C+(L)
@>>> 0$$
where $ I_{\Gamma}$ is the saturated ideal of a curve $\Gamma$ contained
in $C$; as a set, $C$ is the union of $C'$ and $\Gamma$, and we
 have
$$\deg C'+\deg \Gamma =4,$$ and $$g(C')+g(\Gamma)\geq 1-\deg \Gamma.$$
If $\Gamma$ is a line (or equivalently, if $C\cap H_L$ contains a curve of
degree 3), then
$I_{C'}=I_C+(L)$,
$C'$ is a degeneration
 of
a twisted cubic in $H_L$, and $C$ is in $\bH'$.\endproclaim

\demo{Proof} Let $I_{\Gamma}$ be the kernel of the multiplication by $L$
 from
$S$ to $S/I_C$. It is the saturated ideal of a curve $\Gamma$ contained
 in
$C$. Let $K=I_{C'}/I_C+(L)$. Since $C'$ and $C\cap H_L$ differ
 only at a
finite set of points, the sheaf $\widetilde K$ associated to $K$
 has
finite support, and we have, for all $m\in \bN$, the equality:
$$\chi \cO_C(m)=\chi \cO_{\Gamma}(m-1)+\chi \cO_{C'}(m)+lg \ \widetilde K$$
from which we get the relations between the genera and degrees.

Suppose that $\Gamma$ is a line. Then we get $\deg C'=3$ and $g(C')\geq 0$.
(Conversely, if $C\cap H_L$ contains a curve of
degree 3, this curve is equal to $C'$ --- if not, $\deg C'=4$, and $C=C'$ is
contained in a hyperplane --- and $\Gamma$ is a line.)

The ideal $ I_{\Gamma}$ is generated
 by three
independent linear forms $L_1$, $L_2$, $L_3$, so $I_C$ contains
$(LL_1, LL_2, LL_3)$. Moreover, since $I_C$ is generated by six
independent quadratic forms, we have necessarily (for simplicity
 assume
$L=X_4$):
$$I_C=(X_4L_1,X_4L_2,X_4L_3,Q_1,Q_2,Q_3)$$
where $Q_1,Q_2,Q_3$ are three quadratic forms contained in the ideal
$I_{\Gamma}=(L_1,L_2,L_3)$.

If $g(C')=1$, $C'$ is a plane cubic, and the $Q_i$, modulo $X_4$,
 have a
common linear factor $M$. Hence $I_C$ is contained in the ideal $(X_4,M)$,
which is impossible.

It follows that $g(C')=0$, so that $C'$ is a curve of degree 3
 and genus 0
in $H_L\simeq
\bP^3$. Since such a $C'$ is ACM, its ideal $I_{C'}$ in
$k[X_0, \ldots, X_3]$ is generated by three
quadratic forms, which have to be the images of $Q_1,Q_2,Q_3$. Therefore
$I_{C'}=I_C+(X_4)=(Q_1,Q_2,Q_3,X_4)$.

1) Suppose that $\Gamma$ is not in $H_L$. Then $\Gamma$ and $H_L$
intersect in a point
$P$ and the rank of
$L_1,L_2,L_3,X_4$ is 4. In this case $C$ is the scheme theoretic
union of the two curves $C'$ and $\Gamma$ that intersect in $P$ (cf. Lemma 2.1).
 We shall show that $C$ can be deformed into an ACM reduced
curve through $P$.

The set of  curves of degree 3 and genus 0 in
$H_L$ passing through
$P$ is irreducible  of dimension 9 (it is clear if we see it as a
 quotient
of a space of matrices), and it contains smooth twisted cubics. One
can deduce from ([E], Prop.1, p.424), that there
exist, in the  ring
$S[t]$, three polynomials
$ Q_{1,t},  Q_{2,t}, Q_{3,t}$, which are quadratic forms in $X_0,
\ldots, X_3$, and which define a flat family $\cC'$ of curves
contained in $H_L$, parametrized by an affine open subscheme $U$
 of the line,
containing the point $t=0$, such that

-- for $t=0$, $Q_{i,0}-Q_i\in (X_4)$ (and $\cC'_0=C'$),

-- for $t\neq 0$, $\cC'_t$ is a smooth twisted cubic  passing
 through
$P$.

Since $\cC'_t$ goes through $P$, we have $Q_{i,t}\in
(X_4,L_1,L_2,L_3)$.  Write $Q_{i,t}=Q'_{i,t}+X_4 Q''_{i,t}$, where
$Q'_{i,t}\in (L_1,L_2,L_3)$.

Consider the following ideal in $S[t]$:
$$I_t=
(X_4L_1,X_4L_2,X_4L_3,Q'_{1,t},Q'_{2,t},Q'_{3,t}).$$
By Lemma 2.1, for every $t$, $$I_t
=(L_1,L_2,L_3)\cap
(X_4,Q'_{1,t},Q'_{2,t},Q'_{3,t}).$$
This ideal defines a a flat family $\cC$ of curves in $\bP^4$,
 parametrized by  $U$, such that

-- for $t=0$, $\cC_0$ is the curve defined by
$$I_0=
(X_4L_1,X_4L_2,X_4L_3,Q'_{1,0},Q'_{2,0},Q'_{3,0})$$
and since  $Q_i-Q'_{i,0}\in (X_4)\cap
(L_1,L_2,L_3)=(X_4L_1,X_4L_2,X_4L_3)$, we have $\cC_0=C$.

-- for $t\neq 0$, $\cC_t$  is ACM (since the limit $\cC_0=C$ is
ACM) and is equal to the union of the line defined by
$(L_1,L_2,L_3)$ and  the smooth twisted cubic $\cC'_t$,
which intersect at $P$. Hence $\cC_t$ is both ACM and reduced and therefore
belongs
to $\bH'$.

2) Suppose that $\Gamma$ is  in $H_L$, i.e., the rank of $L_1,L_2,L_3,X_4$
is 3. Let $L'$ be a linear form independent of $L_1,L_2,L_3,X_4$,
 and
suppose that
$L_1,L_2,X_4$ are also independent. Since $C$ contains $\Gamma$, we have
$Q_{i}\in (X_4,L_1,L_2)$.  Write $Q_{i}=Q'_i+X_4 L''_{i}$,
where $Q'_{i}\in (L_1,L_2)$.

Consider the following ideal in $S[t]$:
$$I_t=
(X_4L_1,X_4L_2,X_4(X_4+tL'),Q_1+tL'L''_1,Q_2+tL'L''_2,Q_3+tL'L''_3).$$
This ideal defines a flat family $\cC$ of curves in $\bP^4$,
parametrized by  the affine line, such that

-- for $t=0$, $\cC_0=C$

-- for $t\neq 0$, $\cC_t$ is ACM (since the limit $\cC_0=C$
 is
ACM) and is of the same type as the limit curve of the family
 studied
in the preceding case, since the rank of
$L_1,L_2,X_4+tL',X_4$ is 4. Therefore,  $\cC_t$ belongs to
 $\bH'$, and hence the limit curve $C$ does too.\qed \enddemo

\proclaim{Corollary 2.3} If $C$ is  an ACM curve of $\bH$, which is
the union of a double line and a reduced
curve of degree 2, then $C$ is in $\bH'$.\endproclaim

\demo{Proof} The curve $C_{red}$ is either the connected union of three
 lines or
the  connected union of a line and a plane conic, hence is contained
 in a
hyperplane. Apply Proposition 2.2 to this hyperplane.\qed \enddemo

\vskip 0.3cm
We shall now go to the case II.1.

The idea is to deform the ``double conic'' to the union of two
 reduced
conics intersecting in a point.

\proclaim{Lemma 2.4} Let $L$, $M$, $L'$, $M'$ be four independent linear
forms, $Q\in (L',M')$ and $Q'\in (L,M)$ two quadratic forms. Then we have:
$$(L,M,Q)\cap (L',M',Q')=(L,M)\cdot (L',M')+(Q,Q')$$.\endproclaim

\demo{Proof} Write $Q=AL'+BM'$, $Q'=A'L+B'M$.
Then
$\alpha L'+\beta M'+\gamma( A'L+B'M)\in (L,M,AL'+BM')$ if and only
 if there
exists $\lambda$ such that $(\alpha -\lambda A)L'+(\beta -\lambda
 B)M'\in
(L,M)$. This holds if and only if there
exist $\lambda$ and $\mu$ such that $\alpha -\lambda A-\mu M'$ and
$\beta -\lambda B +\mu L'$ belong to $(L,M)$. In this case, we have
 $$\alpha
L'+\beta M'+\gamma (A'L+B'M)-(\lambda Q+\gamma Q')\in (L,M)\cdot (L',M').
\qed$$\enddemo

\proclaim{Proposition 2.5} If $C$ is an ACM curve of $\bH$ which is
 a double
structure on a reduced conic, then $C$ is in $\bH'$.\endproclaim

\demo{Proof}  Suppose that the ideal of $C_{red}$ is $(X_0,X_1,Q)$, where
 $Q$ is
a quadratic form. Since $I_C$ is generated by quadratic forms, we
 may
assume $Q\in I_C$. Let
$L$ be a linear form in
$X_0,X_1$. By Proposition 2.2 (and with the same notation) we may assume
 that
$\Gamma$ is not a line. Therefore we must have $\Gamma=C_{red}$, and
$LI_{\Gamma}=L(X_0,X_1,Q)\subset I_C$. Since this holds for every
 $L$, we
have that $(X_0,X_1)^2\subset I_C$.

Therefore we can write
$I_C=(X_0^2,X_0X_1,X_1^2,Q,X_0L+X_1M,X_0L'+X_1M')$, where $L,L',M,M'$
 are
linear forms in $X_2,X_3,X_4$, and the six quadratic forms are independent.

If $ML'-M'L=0$, we have either $M'=aM$ and $L'=aL$, with $a\in
 k$, but
then
$X_0L+X_1M$ and $X_0L'+X_1M'$ are not independent; or $M=aL$ and
$M'=aL'$, with $a\in k$, and we get
$$(X_0,X_1,L,L')\subset (I_C :(X_0+aX_1))=(X_0,X_1,Q),$$
hence $L=L'=M=M'=0$.

Hence $ML'-M'L\neq 0$. Since $ML'-M'L\in (I_C :(X_0))=(X_0,X_1,Q)$,
we may write
$Q=AX_0+BX_1+ML'-M'L$.

We shall now see that $C$ can be deformed into an ACM
reduced curve.

If the rank of $L,M$ is 2, consider the following ideal in $S[t]$:
$$I_t=
(X_0,X_1)\cdot(X_0-tM,X_1+tL)+(Q,X_0L+X_1M,L'(X_0-tM)+M'(X_1+tL)).$$

This ideal defines a flat family $\cC$ of curves  contained
 in
$\bP^4$ parametrized by  the affine line, such that

-- for $t=0$, $\cC_0=C$.

-- for $t\ne 0$, $\cC_t$ is ACM and is defined by the ideal
$$\multline
(X_0,X_1)\cdot (X_0-tM,X_1+tL)\\
\qquad +(L'(X_0-tM)+M'(X_1+tL),X_0(L'+tA)+
X_1(M'+tB)).
\endmultline
$$
It follows from Lemma 2.4 that $\cC_t$ is the
union of the two conics defined by the ideals
$(X_0,X_1,L'(X_0-tM)+M'(X_1+tL))$ and
$(X_0-tM,X_1+tL, X_0(L'+tA)+X_1(M'+tB))$. The first conic is $C_{red}$,
and the second is reduced for general $t$, because it varies in a
 flat
family of conics which is a deformation of $C_{red}$.

If the rank of $L',M'$ is 2, the situation is similar.

If the rank of $L,M$ and  the rank of $L',M'$ are both 1, we can
 suppose
that $M=L'=0$ and that $C$ is a double
structure on the degenerated conic defined by $(X_0,X_1,LM')$. Choose
 a
linear form $L''$, independent of $M'$, and  consider the following
 ideal
in
$S[t]$:
 $$I_t=(X_0^2,X_0X_1,X_1^2,AX_0+BX_1-M'L,X_0L,tX_0L''+X_1M')$$

It defines a a flat family $\cC$ of curves  contained in $\bP^4$
parametrized by  the affine line, such that

-- for $t=0$, $\cC_0=C$,

-- for $t\ne 0$, $\cC_t$ is  ACM, and it is of the same type
(with $L,M,L',M'$ replaced by $L,0,tL'',M'$)
as the limit curve of the family studied in the preceding case, since
 the
rank of
$M',L''$ is 2. Hence, $\cC_t$ belongs to
 $\bH'$, and so does the limit curve $C$.\qed \enddemo

\vskip 0.3cm

For the remaining cases, II.2 (union of a reduced  line
and a triple structure on another line) and III (quadruple structure
on a line), we need another approach, using projections.
\vskip 0.3cm
\vskip 0.3cm

\head{3. Projection on a hyperplane}\endhead
Let $P$ be a point of ${\bP^4}$ and $\pi$ the projection  on a hyperplane
$\simeq {\bP^3}$ from the point $P$. Let  $U={\bP^4} - \{P\}$ denote the
domain of
definition of $\pi$. If $X$ is a closed subscheme of ${\bP^4}$
not passing through $P$, the scheme theoretic image $X'$
of $X$ is well defined: the ideal sheaf
$\cI_{X'}$ is the
kernel of the composed map $\cO_{\bP^3}\to   \pi_* \cO_U
\to \pi_* \cO_X$. The  closed subscheme $X'$ of ${\bP^3}$ is
the projection of $X$ and there is a morphism from $X$ to $X'$
induced by $\pi$. Moreover, the projection of $X_{red}$ is $X'_{red}$.

\proclaim{Lemma 3.1} Let $C$ be a (locally
Cohen-Macaulay) curve in ${\bP^4}$, not contained in a hyperplane,
and $\pi$ a projection  on a hyperplane from a point not on $C$.
 Then the projection $C'$ of $C$ is a (locally Cohen-Macaulay) curve not
contained in a plane.\endproclaim

\demo{Proof} First, $C'$ is not plane because $C$ is not contained
in a hyperplane.

There is a natural injection $ \cO_{C'} \to
\pi_* \cO_C$. If  $C'$ is not locally Cohen-Macaulay, $C'$ has
an associated component of dimension 0, therefore $C$ has a component
which is a line
contracted by $\pi$, and the center of the projection is on $C$,
 hence we get a contradiction.\qed \enddemo

\proclaim{Proposition 3.2} Let $C$ be an ACM curve of $\bH$, and let
$\pi$ be a
projection with center $P\notin C$.
 The projection $C'$ of $C$ is a curve, not contained in a plane, of genus
0 or 1, and
degree 3 or 4. If $C_{red}$ is a union of lines, then
$\pi:C\rightarrow C'$ is an isomorphism at every point $x\in C$
 such that
$P\notin T_{C,x}$.\endproclaim

\demo{Proof}  From the exact sequence
$$0\to \cO_{C'} \to
\pi_* \cO_C \to \cK \to 0$$
 we obtain the following equality for all $m\in \bZ$
$$\chi  \cO_{C'}(m) + \chi  \cK(m)=\chi \pi_*
\cO_C(m)=\chi  \cO_C(m),$$
because $\pi_{\vert C}$ is finite. Therefore the  degree $d'$ of $C'$
is $\leq 4$, and its genus $g'$ satisfies $$g'\leq {{(d'-2)(d'-3)}\over 2}.$$

Now we have $\dim H^0(C',\cO_{C'})\leq
\dim H^0(C', \pi_* \cO_C)=\dim H^0(C, \cO_{C})=1$.
So we get $\dim H^0(C', \cO_{C'})=1$, hence $g'=\dim
H^1(C', \cO_{C'})\geq 0$. Since $C'$ is not plane,
we have $d'\geq 3$.

Suppose $C_{red}$ is a union of lines. Since the projection of a line is a
line,
$\pi$ is set theoretically bijective, therefore it is an isomorphism at a point
$x$ if and only if it is not ramified at $x$, that is, if and only if
$P\notin T_{C,x}$.\qed \enddemo

\proclaim{Proposition 3.3}. Let $C$ be an ACM curve of $\bH$, and let
$\pi$ be a
projection with center $P\notin C$.
 If the projection $C'$ of $C$ is a curve of
degree 4 and genus 0, then the curve $C$ belongs to $\bH'$.\endproclaim

\demo{Proof} Suppose that $C'$ is a (4,0). Then for all $m\in \bZ$, we have
 $\chi
(\cK(m))=0$, hence $\cK=0$ and $\cO_{C'}\simeq
 \pi_* \cO_C$, so the restriction of $\pi$ to $C$ is an embedding. The
curve $C'$
is not linearly normal in $\bP^3$: the sheaf $\cO_{C'}(1)$ has five
sections which correspond to the embedding of $C'$ as $C$ in $\bP^4$.

Now, we know from [M-D--P] that $C'$ is a curve
in the biliaison class of two disjoint lines: in particular,
$h^1(\cI_{C'}(m))=0$ for $m\ne 1$, and $h^1(\cI_{C'}(1))=1$.
It follows, on the one hand, that $h^0(\cI_{C'}(2))=1$, so that $C'$
is contained in a quadric, and on the other hand that
$h^1(\cI_{C'}(2))=h^2(\cI_{C'}(1))=0$, hence --- by the
Castelnuovo--Mumford regularity criterion --- that $C'$ is the intersection of
cubic
surfaces. In particular, one can make a linkage of type (2,3), and one obtains
a curve of degree $2$ and genus $-1$ which is either the union of two disjoint
lines or a double structure on a line.\enddemo

To conclude we need the  following result.

\proclaim{Lemma 3.4} Let $C'$ be a curve of degree $4$ and genus $0$ in $\bP^3$
such that $C'_{red}$ is a union of lines. Suppose $C'$ is contained in
a complete intersection of type (2,3). Then there exists a curve $C'_1$
of degree $3$ and genus $0$ and a line $D'$, both contained in $C'$, and
an exact sequence
$$0\to \cO_{D'}(-1) \to
\cO_{C'} \to \cO_{C'_1} \to 0.$$\endproclaim

Let us assume this result. Because $C'$ can be embedded in $\bP^4$, so can
$C'_1$, and its image is a curve $C_1$ of degree $3$ contained in $C$.
Moreover, the image of $\cO_{D'}(-1) \to \cO_{C'}$ corresponds to a
section of $\cO_{C}(1)$ vanishing on $C_1$, hence it comes from a linear
form defining a hyperplane $H$ in $\bP^4$ such that $H\cap C$ contains $C_1$.
We can therefore apply Proposition 2.2.

\medskip
Lemma 3.4 follows from another lemma:

\proclaim{Lemma 3.5} Let $\Gamma$ be a curve of degree $2$ and genus $-1$ in
$\bP^3$. Assume that $\Gamma$ can be linked by a complete intersection $X$ of
type (2,3) to a curve supported by lines. Then there exists a curve
$\Gamma_1$ of degree $3$ and genus $0$ containing $\Gamma$ and contained in
$X$.\endproclaim

\demo{Proof} Let $X_0, X_1, X_2, X_3$ denote coordinates of $\bP^3$. Recall (cf.
[Mi]) that, up to projective equivalence, the equations of a structure of degree
$2$ and genus
$-1$  on a line are $(X_0^2, X_0X_1, X_1^2, X_0X_2+X_1X_3)$. One can then
verify that the singular quadric surfaces which contain such a curve are
unions of two planes or double planes.

We shall have to distinguish between
several cases. Let
$Q$ (resp.
$S$) denote the quadric (resp. cubic) surface (and also its equation), and set
$X=Q\cap S$.

1) If $Q$ is smooth, $\Gamma$ is linearly equivalent to two lines of one system
on $Q$, and $X$ contains at least one line $D$ from the other system --- hence
$\Gamma_1 = \Gamma\cup D$ works.

2) If $Q$ is a cone, it does not contain two disjoint lines, nor a double line
of genus $-1$.

3) If $Q$ is the union of two planes $H$ and $H'$, and $\Gamma$ is the union
of two lines $D$ and $D'$, we may assume $D\subset H$ and $D'\subset H'$. Set
$\Delta = H\cap H'$. If $X$ contains $\Delta$, $\Gamma_1 = \Gamma\cup
\Delta$ works. If not, we have $S\cap H = D\cup D_1 \cup D_2$ and
$S\cap H' = D'\cup D'_1 \cup D'_2$, where $D_1, D_2, D'_1, D'_2$ are (not
necessarily distinct) lines. Then $D'\cap \Delta$ is contained in
$(D\cup D_1 \cup D_2)\cap \Delta$, hence for example in $D_1\cap \Delta$.
So $\Gamma_1= \Gamma \cup D_1$ works.

4) Suppose $Q$ is the union of two planes $H$ and $H'$ and $\Gamma$ is
a double line. We may assume  that
$$I_\Gamma = (X_0^2, X_0X_1, X_1^2, X_0X_2+X_1X_3)$$ and $Q=X_0X_1$.
Then $$S\in (X_0, X_1LL')\cap (X_1, X_0MM'),$$ where $L, L'$ (resp.
$M, M'$) are two linear forms independent of $X_0$ (resp. $X_1$), and hence
$$S=AX_0X_1+\alpha X_0MM'+\beta X_1LL',$$ where $A$ is a linear form and
$\alpha\beta\ne 0$. Since $S\in I_\Gamma$, we must have $X_1LL'\in
(X_1^2, X_3)$, hence $LL'\in (X_1, X_3)$ and, say, $L\in (X_1, X_3)$.

If $L\notin (X_1)$, let $D$ denote the line defined by $(X_0, L)$ and
set $\Gamma_1=\Gamma \cup D$. Then $$I_{\Gamma\cap D}=I_\Gamma + (X_0, L)=
(X_0, L, X_1^2, X_1X_3)=(X_0, L, X_1^2),$$ hence $\Gamma_1$ is a
curve of degree $3$ and genus $0$.

If $L\in (X_1)$ and $M\in (X_0)$, then $S\in (X_0^2, X_0X_1, X_1^2)$ and this
ideal is the ideal of a curve $\Gamma_1$ that works.

5) If $Q$ is a double plane, then $\Gamma$ is a double line. We may assume
$Q=X_0^2$ and $$I_\Gamma=(X_0^2, X_0X_1,X_1^2,X_0X_2+X_1X_3).$$
As in 4) we get $S\in (X_0, X_1LL')$ and $LL'\in (X_1, X_3)$ and, say,
$L\in (X_1, X_3)$.

If $L\notin (X_1)$, we conclude as in 4). If $L\in (X_1)$, then
$$S\in (X_0, X_1^2)\cap I_\Gamma = (X_0^2, X_0X_1, X_1^2)$$ and again we
conclude
as in 4).\qed \enddemo

\medbreak

\demo{Proof of 3.4} By assumption, $C'$ is linked to a curve $\Gamma$
via a complete intersection $X$ of type (2,3), so $\Gamma$ has degree $2$
and genus $-1$.
With the notations as in Lemma 3.5,
the curve $\Gamma_1$ is linked to a curve $C_1'$ contained in $C'$.
Because $\Gamma_1$ is of degree $3$ and genus $0$, so is $C_1'$.

Set $\cF = \cI_{C'_1}/\cI_{C'}$. For all $m$, we have
$\chi (\cF(m))=m$. Moreover, $h^0(\cF)=0$, $h^0(\cF(1))=1$,
$h^1(\cF)=0$, so that $\cF(1)$ is generated by its sections and
$\cF \simeq \cO_{D'}(-1)$, for some $D'$ contained in $C$.
It is then easy to see that $D'$ is a line.\qed \enddemo

\vskip 0.3cm
It remains to study the curves $C$ in the cases II.2 and III
($C_{red}$ is then a union of lines) such that for {\it every} projection
$\pi$ from a point not on $C$, the scheme theoretic image $C'$ of $C$
is a (3,0) or a (4,1). We may also assume (cf. Proposition 1.3) that
the tangent spaces $T_{C,x}$ at the points $x$ of a component of $C$
are not all of dimension 4. We deduce from Propositions 3.2 and 3.3 the
following result :

\proclaim{Corollary 3.6} Let $C$ be an ACM curve of $\bH$ such that
$C_{red}$ is a union of lines, and such that for every projection
$\pi$ from a point $P\notin C$, the scheme theoretic image $C'$ of
$C$ is a (3,0) or a (4,1). Then the union of the tangent
spaces $T_{C,x}$ at all points $x$ of $C$ is the whole space. Moreover, if
$P\in \cap \ T_{C,x}\backslash C$, then $C'$ is a (3,0); if $P\notin\
\cap T_{C,x}$, then there exists
 a unique
point $x\in C$ such that $P\in T_{C,x}$, $\pi$ is an isomorphism
outside of $x$, and $C'$ is a (4,1).\endproclaim

\demo{Proof} In the exact sequence
$$0\to \cO_{C'} \to
\pi_* \cO_C \to \cK \to 0$$
the support of $\cK$ is the set of points of $C$ where $\pi:
C\to C'$ is not an isomorphism. Proposition 3.2 shows that this is also
the set of points $x$ of $C$ such that $P\in T_{C,x}$. If this set is
empty,
$C'$ is a (4,0).

The degree of $C'$ is 4 if and only if the support
of $\cK$ is finite, and this is
equivalent to $P\notin \cap T_{C,x}$.
In this case, for all $m\in \bZ$, we have
 $\chi (\cK(m))=1$ and the support of $\cK$ is one point.\qed \enddemo

\medskip

We need now a description of the curves of type (4,1) or (3,0) in $ \bP^3$,
supported by
lines. The computations are elementary and will not be given.

\proclaim{Lemma 3.7} Let
$\Gamma\subset\bP^3$ be a curve of type (4,1) which  is the union of a reduced
line $X_1=X_3=0$ and a triple line with support $X_1=X_2=0$. Then the ideal of
$\Gamma$, up to a linear change of coordinates, is  given as
$$(X_1X_2, X_1^2+X_2X_3) \quad \text{or}  \quad (X_1X_2,X_2X_3,X_1^3)
\quad \text{or}  \quad (X_1^2, X_1X_2,X_1Q+X_2^2X_3)$$ where $Q$ is a
quadratic form independent of $X_1,X_2$.
In the first two  cases, the tangent spaces to $\Gamma$ are constant along the
triple line, equal to
a plane.
\indent Moreover, in the first case, $\Gamma$ links the reduced line to
the triple line given by the ideal $(X_1^2+X_2X_3,X_1X_2,X_2^2)$, which is a
curve of type (3,0). In the third case, the tangent spaces to $\Gamma$
are also a constant plane, except at the points where $Q$ vanishes. If
$Q=0$, the tangent spaces are of
dimension $3$ (equal to the whole space) along the triple line.\endproclaim

\proclaim{Lemma 3.8} Let $\Gamma\subset\bP^3$ be a curve of type (4,1)
which  is  a
quadruple structure on the line $X_1=X_2=0$, then the ideal of
$\Gamma$, up to a linear change of coordinates, is  given as
$$(X_1^2,X_2^2) \quad \text{or}  \quad (X_1^2,X_2^2+X_1X_3) \quad \text{or}
\quad (X_1^2,X_1X_2,X_1Q+X_2^3)$$
where $Q$ is a quadratic form.
In the second and third cases (if $Q\ne 0$), (almost all) the tangent
spaces are
 constant,
equal to a plane. In the first and third cases (if $Q=0$),
the tangent spaces are  the whole space.
\indent Moreover, in the first (resp. second) case, $\Gamma$ links the
reduced line
 $X_1=X_2=0$
to a triple structure of type (3,0) on the same line, with equations
$(X_1^2,X_1X_2,X_2^2)$ (resp.$(X_1^2,X_1X_2,X_2^2+X_1X_3)$).\endproclaim

As an easy consequence of these two results, we get :

\proclaim{Lemma 3.9} In $\bP^3$, the only (3,0) structures on a line are
(up to a change of
coordinates) the ones given by the ideal $(X_1^2,X_1X_2,X_2^2+\alpha X_1X_3)$
($\alpha =0$ gives the $3$-fold line).\endproclaim

\demo{Proof} A curve $\Gamma$ of type (3,0) is the intersection of quadric
surfaces,
hence is contained
in a curve $\Gamma'$ of type (4,1), which is the complete intersection of
two quadrics. If
$\Gamma$  is supported by a line $D$, $\Gamma'$ is either  a quadruple
structure on $D$ or the union of a triple structure on $D$ with
another reduced line $D'$. In the first (resp. the second)  case, $\Gamma'$
links
$D$ (resp. $D'$) to $\Gamma$ and the equations of $\Gamma$ are given by
Lemmas 3.7 and 3.8.\qed \enddemo

\medskip
\remark{Remark 3.10} More generally, consider the $n$-uple structure on the
line $X_1=\cdots =
X_{n-1}=0$ in $\bP^n$ given by the ideal generated by the $2\times 2$-minors
of the matrix
$$\pmatrix 0&X_1&\ldots &X_{n-2}&X_{n-1}\\
X_1&X_2&\ldots &X_{n-1}&\alpha X_n
\endpmatrix
$$
This curve is ACM (and equal to the $n$-fold line if $\alpha =0$).
For $n=4$, we computed (using the computer program Macaulay) the tangent
space to the
Hilbert scheme at such a point and found its dimension to be $24$ when
$\alpha \neq 0$. Therefore, these curves correspond to singular points on the
Hilbert scheme if $n=4$, and we conjecture the same is true for all $n\geq 4$.
\endremark

\proclaim{Proposition 3.11} Let $D\subset \bP^4$ be a line and let $C$ be
an ACM curve of $\bH$
which is either the union of a triple structure on the line $D$
 and a reduced
line, or a quadruple structure on $D$, and such that for every projection
$\pi$ from a point $P$ not on $C$, the scheme theoretic image $C'$ of
$C$ is a (3,0) or a (4,1). If dim $T_{C,x}=3$ for a
general point $x$ of $D$, then $C$ is in $\bH'$.\endproclaim

\demo{Proof} Assume $D$ is the line $X_1=X_2=X_3=0$. Let $x_1$ and $x_2$ be two
distinct points on $D$. Since
dim $T_{C,x_1}\cap T_{C,x_2}\geq 2$, one may choose a projection
 center
$P\in T_{C,x_1}\cap T_{C,x_2}\backslash C$. It follows from Corollary 3.6 that
$P\in T_{C,x}$ for every point $x$ of $D$,
hence the tangent spaces along $D$ contain a fixed plane, say
$X_1=X_2=0$.

Let $F\in I_C$, $F=A_1X_1+A_2X_2+A_3X_3$. The vector
$\pmatrix A_1&A_2&A_3&0\endpmatrix$
satisfies $A_3=0$ along the line $D$, hence $A_3\in (X_1,X_2,X_3)$,
 so
that $I_C\subset (X_1,X_2,X_3^2)$.

We know that $I_C$ is generated by
six quadrics $Q_1,...,Q_6$. Set $$Q_i= X_1L_i+X_2M_i+\lambda_iX_3^2.$$
Because the general tangent spaces to $C$ are of dimension $3$, the
matrix
$$\pmatrix L_1&\ldots &L_6\\
M_2&\ldots &M_6
\endpmatrix
$$
is of rank $\leq 1$ along the line $D$. The proof of the
following lemma is easy:\enddemo

\proclaim{Lemma 3.12} In the above situation, either
\roster
\item"{(i)}" $L_i=cM_i$ mod $(X_1,X_2,X_3)$ for some constant $c$ and for
all $i$,
\endroster
or
\roster
\item"{(ii)}" there exists $i_0$ and constants $\lambda_i$ such that
$L_i=\lambda_iL_{i_0}$ and $M_i=\lambda_iM_{i_0}$, mod $(X_1,X_2,X_3)$.
\endroster
\endproclaim

\medskip
Let us now return to the proof of Proposition 3.11.
It follows from the lemma that we have two possibilities:

\noindent (i) $L_i=cM_i$ mod $(X_1,X_2,X_3)$.

Then $Q_i=M_i(cX_1+X_2)$  mod $(X_1,X_2,X_3)^2$, so
$$I_C+(cX_1+X_2)\subset (cX_1+X_2)+(X_1,X_2,X_3)^2.$$
The hyperplane section of $C$ given by $cX_1+X_2=0$ contains a
 curve of degree
$3$, and we can therefore conclude using Proposition 2.2.
\medskip
\noindent (ii) We may assume $(L_i,M_i)\subset (X_1,X_2,X_3)$ for
 $i\ne 1$,
hence $I_C\subset (Q_1)+(X_1,X_2,X_3)^2$, with $Q_1=X_1L_1+X_2M_1$,
and where $L_1$ and $M_1$ are independent of $(X_1,X_2,X_3)$.

If $M_1=0$ (or $L_1=0$), or, more generally, if $M_1=\lambda
 L_1$, then
$$I_C+(X_1)\subset (X_1)+(X_2,X_3)^2,$$ and we may conclude as above
 by Proposition 2.2.

Suppose $L_1$ and $M_1$ are independent. We can suppose $L_1=X_0$ and
$M_1=X_4$.

If $C$ contains a reduced line $D'$, which meets $D$ in a
point $(\alpha_0,0,0,0,\alpha_4)$, $D'$ is contained in the tangent
space to $C$ at this point, which is defined by $\alpha_0
X_1+\alpha_4X_2=0$. Then

$$I_C+(\alpha_0 X_1+\alpha_4 X_2)\subset (X_1,X_2,X_3^2)\cup I_{D'},$$
and we may conclude as above by Proposition 2.2.

If $C$ is a quadruple structure on the line $D$, project
$C$ from the point
$(0,0,0,1,0)$ (which belongs to all the tangent spaces to $C$)  into
the plane $X_3=0$. We obtain a curve of type (3,0) supported on
the line
$X_1=X_2=0$, such that its ideal is generated by three quadratic
 forms and
is contained in $$((Q_1)+(X_1,X_2,X_3)^2)\cap k[X_0,X_1,X_2,X_4]
=(Q_1)+(X_1,X_2)^2.$$
>From what we have seen concerning curves of type (3,0) supported
on a
line, and because $Q_1$ is irreducible, the image curve must be the
$3$-fold line. Hence $(X_1,X_2)^2\subset I_C$.

If instead we project from the point $(0,0,1,0,0)$ (which does not
belong to all the tangent spaces)  into the plane $X_2=0$, the
image curve is of type (4,1), supported on a line, and the general tangent
space has dimension 3:

-- if its ideal has the form $(X_1^2, X_1 L, L^3)$, where $L$ is a
linear form in $X_1, X_3$, linearly independent of $X_1$, we see that
$(I_C:X_1)$ contains $X_1,X_2,L$ and we apply Proposition 2.2, intersecting
$C$ with the hyperplane $ X_1=0$.

-- if its ideal is generated by the squares of two linear
forms, we can suppose that the intersection
$I_C\cap k[X_0,X_1,X_3,X_4]$, which already contains $X_1^2$, also
contains $X_3^2$.  Then $(X_1,X_2)^2+(X_3)^2\subset I_C\subset
(X_1,X_2,X_3^2)$. The two ``missing'' generators are $Q_1+Q'_1$ and $Q'_2$,
with $(Q'_1,Q'_2)\subset (X_1,X_2,X_3)^2$. We may therefore choose
$Q'_2=\alpha X_1X_3+\beta X_2X_3$ and apply Proposition 2.2, intersecting
$C$ with the hyperplane $\alpha X_1+\beta X_2=0$.\qed

\proclaim{Proposition 3.13} Let $D\subset \bP^4$ be a line, and let
$C$ be  an ACM curve of $\bH$
which is either the union of a triple structure on the line $D$ and a reduced
line, or a quadruple structure on $D$, and such that for every projection
$\pi$ from a point $P\notin C$, the scheme theoretic image $C'$ of
$C$ is a (3,0) or a (4,1). If dim $T_{C,x}=2$ for a
general point $x$ of $D$, then $C$ is in $\bH'$.\endproclaim

\demo{Proof}  Assume $D$ is the line $X_1=X_2=X_3=0$. A general projection of
$C$ is of type (4,1), with $2$-dimensional
tangent spaces along the multiple line. By Corollary 3.8, these tangent
planes are
equal.

If the tangent planes to $C$ are not constant, the hyperplane they span
must contain the center of projection (since they all project to
the same plane) --- but this cannot happen for a {\it general}
projection. Therefore they are constant, equal to, say, $X_1=X_2=0$.

As in the proof of Proposition 3.11, one has $I_C\subset (X_1,X_2,X_3^2)$.
We
can therefore write the six quadratic generators of $I_C$ on the form
$Q_i=X_1L_i+X_2M_i+Q'_i$, where $Q'_i\in (X_1,X_2,X_3^2)$ and
$L_i,M_i\in (X_0,X_4)$.

If the rank of $(L_1,...,M_6)$ is equal to $2$, then we see from the equations
that there is no point $x$ on $C$ where the tangent space has dimension 4, and
only a finite number of points where the tangent space has dimension 3, so the
union of the tangent spaces to $C$ along the multiple line is different from
$\bP^4$, and we get a contradiction.
So we may assume $Q_i=X_0(\alpha_iX_1+\beta_iX_2)+Q'_i$.

If the rank of $(Q'_1,...,Q'_6)$ is $\leq 5$, we may assume that
for some $i$, we have $Q'_i=0$, hence $X_0(\alpha_iX_1+\beta_iX_2)
\in I_C$. Then $X_0$ is a zero divisor in $S/I_C$, and necessarily,
$C$ is the union of a triple line $A$ and a reduced line $D$
contained in $X_0=0$. Moreover, $\alpha_iX_1+\beta_iX_2 \in I_A$,
 hence
$I_C+(\alpha_i X_1+\beta_i X_2)\subset I_A$. When we intersect $C$
 by
$\alpha_i X_1+\beta_i X_2=0$, we get a curve of degree $\geq 3$,
 hence again Proposition 2.2 applies.

If not, we may assume that $I_C$ is generated by $Q_{ij}=
X_0L_{ij}+X_iX_j$, where $L_{ij}=\alpha_{ij}X_1+\beta_{ij}X_2$, for
$i,j\in \{1,2,3\}$. Hence $I_C$ contains
$X_0(X_iL_{ij}-X_jL_{ii})$. Since $I_C+(X_0)=(X_0)+(X_1,X_2,X_3)^2$, $X_0$
is not a zero divisor in $S/I_C$. Hence $I_C$ contains $X_iL_{ij}-
X_jL_{ii}$.

If the rank of $L_{11},L_{12},L_{22}$ is $\leq 1$, then $I_C$ contains
two quadratic forms in $X_1,X_2$, and hence the square of some linear form,
and we may assume, e.g., that $L_{11}=0$. Then
$(X_1L_{12},X_1L_{13},X_1^2)=(X_1^2,\beta_{12}X_1X_2,\beta_{13}X_1X_3)
\subset I_C$. Hence $ I_C$ contains $X_1(X_2+\lambda_2 X_0)$ where $\lambda_2=0$
if $\beta_{12}\ne 0$, and $\lambda_2=\alpha_{12}$ if $\beta_{12}=0$. For the
same reason,  $ I_C$ contains $X_1(X_3+\lambda_3 X_0)$. So we may intersect with
$X_1=0$ and apply Proposition 2.2.

If the rank of $L_{11},L_{12},L_{22}$ is $2$, then $I_C$ contains
a quadratic form in $X_1,X_2$. If this is a square, we conclude as
in the preceding case. If it is not a square, we may assume
$L_{12}=0$, hence $(L_{11}X_2,L_{22}X_1,X_1X_2)=(\beta_{11}X_2^2,
\alpha_{22}X_1^2,X_1X_2)\subset I_C$. But because $I_C$ does not
 contain
a square in $X_1,X_2$, we have $\beta_{11}=\alpha_{22}=0$. Hence
$$Q_{11}=X_1(X_1+\alpha_{11}X_0), \quad Q_{22}=X_2(X_2+\beta_{22}X_0),$$
so that $X_1+\alpha_{11}X_0$ and $X_2+\beta_{22}X_0$ are
zero divisors in $S/I_C$ (note that $\alpha_{11}\beta_{22}\neq 0$ since the rank
of $L_{11},L_{12},L_{22}$ is $2$). It follows that
$C$ must contain a reduced line $D'$,
and that $(X_1 X_2,X_1+\alpha_{11}X_0,X_2+\beta_{22}X_0)\subset I_D'$.
But this implies $X_0\in I_D$, contradicting the fact that $X_0$
is not a zero divisor in $S/I_C$.\qed \enddemo
\medskip

We have now proved the following.
\medskip
\proclaim{Proposition 3.14} If $C$ is an ACM curve of $\bH$ and $C$ is
equal either to the union of a reduced  line
and a triple structure on another line, or to
a quadruple structure on a line, then $C$ is in $\bH'$.\endproclaim

\medskip
Putting Proposition 3.14 together with Corollary 2.3 and
Proposition 2.5, we obtain the following
result.
\medskip

\proclaim{Theorem 3.15} The ACM curves form an irreducible open subscheme
of the Hilbert scheme $Hilb^{4n+1}(\bP^4)$.\endproclaim

As noted in Section 1, the $4$-fold lines correspond to points on
the 21-dimensional scheme $\bH^0$
where the tangent space has dimension 36. Since these points are not
contained in any other component of the Hilbert scheme, they must be
singular points on $\bH^0$. Hence the singular locus of $\bH^0$
contains a variety isomorphic to Grass(1,4). In fact,
in Remark 3.10 we described more general ACM $n$-fold structures
on a line; for
$n=4$ we checked that these curves also correspond
to singular points on the Hilbert scheme. Hence we obtain the following
corollary.

\medskip
\proclaim{Corollary 3.16} The irreducible component $\bH'$ of
$Hilb^{4n+1}(\bP^4)$ containing the rational normal curves is singular. More
precisely, the open subscheme $\bH^0\subset \bH'$ corresponding to
ACM curves
has a singular locus that contains (strictly) a variety isomorphic to
Grass(1,4).
\endproclaim

\head{4. ACM curves in $\bP^4$ with Hilbert polynomial $4m+1$}
\endhead
The study made in the previous section allows us to give a precise description
of the curves corresponding to points of $\bH^0$.

\medskip
\proclaim{Proposition 4.1} Let $C$ be a curve of $\bH^0$. Then one of the
following holds:
\roster
\item"{(i)}" $C$ is a rational normal curve.
\item"{(ii)}" $C$ is the union of two smooth conics intersecting in one point.
\item"{(iii)}" $C$ is a double structure on a conic.
\item"{(iv)}" There exists a hyperplane $H$ such that $C\cap H$ contains a
curve of degree $3$.
\endroster
\endproclaim

\demo{Proof} If $C$ is reduced and irreducible, $C$ is a RN curve. If $C$
is reduced,
but is neither irreducible nor the union of two smooth conics, $C$ is the union
of a line $D$ and a connected reduced (possibly reducible)
curve $C'$ of degree $3$ that intersects $D$. Hence we're in the case (iv).

Assume $C$ is not reduced and not equal to a double structure on a conic.
The $C$ is either the union of a reduced conic with a double structure on
a line, or the union of a line and a triple structure on
another line, or a quadruple structure on a line. By the proofs of Corollary
2.3 and Proposition 3.14, (iv) holds in these cases.\qed \enddemo

\medskip

We shall now give the corresponding structures, i.e., describe the ideals
of these curves. The following result --- apart from (i), which is classical
--- is a direct consequence of Sections 2 and 3.

\proclaim{Proposition 4.2} With the notations of Proposition 4.1, the ideal
$I_C$ of the curve $C$ is given, up to projective equivalence, as
\roster
\item"{(i)}" $I_C$ is the ideal generated by the $(2\times 2)$-minors of
the matrix
$$\pmatrix X_0&X_1&X_2&X_3\\
X_1&X_2&X_3&X_4
\endpmatrix
$$
\item"{(ii)}" $I_C=(X_0,X_1)\cdot (X_2,X_3)+(Q,Q')$, where $Q$ and $Q'$ are
quadratic forms such that $Q\in (X_2,X_3)$, $Q\notin (X_0,X_1)$, and
$Q'\in (X_0,X_1)$, $Q'\notin (X_2,X_3)$.
\item"{(iii)}" $I_C=(X_0^2, X_0X_1,X_1^2,X_0L+X_1M,X_0L'+X_1M',
AX_0+BX_1+ML'-M'L)$, where $A,B,L,M,L',M'$ are linear forms in $X_2,X_3,X_4$
such that $ML'-M'L\ne 0$.
\item"{(iv)}" $I_C=(LX_1,LX_2,LX_3,Q_1,Q_2,Q_3)$, where $L$ is
linear and $Q_1,Q_2,Q_3$ are quadratic,
$(Q_1,Q_2,Q_3)\subset (X_1,X_2,X_3)$, and the
ideal $(L,Q_1,Q_2,Q_3)$ defines a $(3,0)$ in the hyperplane $L=0$.
\endroster
\endproclaim

In order to complete the description of the points of $\bH^0$ we must
verify that all the above structures define ACM curves with Hilbert polynomial
$4m+1$.

The cases (i), (ii), (iii) are verified by considering the generators
of the ideal $I_C$ (this can be done directly, or by using the computer
program Macaulay). To treat the last case, we need the following result.

\proclaim{Lemma 4.3} Let $L_1, L_2, L_3$ (resp. $Q_1, Q_2, Q_3$) be independent
linear (resp. quadratic) forms. Assume that $(Q_1, Q_2, Q_3)\subset
(L_1, L_2, L_3)$, and let $L$ be a linear form such that the ideal
$(L,Q_1, Q_2, Q_3)$ defines a
curve $C'$ of degree $3$ and genus $0$. Then the following are equivalent:
\roster
\item"{(a)}" The curve $C$ defined by the ideal
$I_C=(LL_1,LL_2,LL_3,Q_1,Q_2,Q_3)$ is ACM
(of degree $4$ and genus $0$).
\item"{(b)}" The ideal $((Q_1,Q_2,Q_3):L)$ is contained in $(L_1, L_2, L_3)$.
\endroster
\endproclaim

\demo{Proof} Let $D$ denote the line defined by the ideal $I_D=(L_1,L_2,L_3)$.
We note that (b) is equivalent to the equality
$(I_C:L)=I_D$. Hence (a) implies (b) by Proposition 2.2.

Conversely, if $(I_C:L)=I_D$, then there is an exact sequence
$$0 @>>> S/I_D(-1) @> {\cdot L} >> S/I_C @>>> S/I_{C'}
@>>> 0,$$ where $S/I_D$ and $S/I_{C'}$ are Cohen--Macaulay of dimension
$2$ --- hence so is $S/I_C$.\qed \enddemo

\medskip
\remark{Remark 4.4} If $L,L_1, L_2, L_3$ are linearly independent, (b)
holds. In fact, we then have $$((Q_1,Q_2,Q_3):L)\subset ((L_1, L_2, L_3):L)
=(L_1, L_2, L_3).$$
Geometrically, we can state this as folows:.
Let $H\subset \bP^4$ be a hyperplane,  $D$  a line not
contained in $H$, and $C'\subset H$ a curve of degree 3 and genus 0
passing through the point $H\cap D$. Then the scheme theoretic union $C$ of
$D$ and $C'$ is an ACM curve of $\bH$.\endremark
\medskip

To sum up, we have shown that the ACM curves can be described as follows:

\proclaim{Theorem 4.5} The ACM curves of $\bP^4$ with Hilbert polynomial
$4m+1$ are the curves of the following four types:
\roster
\item"{(i)}" a rational normal curve
\item"{(ii)}" the union of two smooth conics intersecting in one point and not
contained in a hyperplane
\item"{(iii)}" a double structure on a conic, as described in
Proposition 4.2 (ii).
\item"{(iv)}" a curve obtained from a curve
of degree $3$ and genus $0$ contained in a hyperplane, and a line, satisfying
the conditions of Lemma 4.3.
\endroster
\endproclaim

\medskip
\Refs
\roster
\item"{[C1]}" J. A. Christophersen, {\it Notes on the component of the
rational normal curves
in $Hilb_{nt+1}$,} manuscript, 1989.

\item"{[C2]}" J. A. Christophersen, pvt. comm.

\item"{[E--R--S]}" D. Eisenbud, O. Riemenschneider, F.-O. Schreyer, {\it
Projective
resolutions of Cohen--Macaulay algebras,} Math. Ann. {\bf 257} (1981), 85--98.

\item"{[E]}" G. Ellingsrud, {\it Sur le sch\'ema de Hilbert des
vari\'et\'es de codimension 2
dans $\bP^e$ \`a c\^one de Cohen--Macaualay,} Ann. scient. \'Ec. Norm. Sup.
{\bf 8}
(1975), 423--432.

\item"{[G]}" A. Grothendieck, {\it Techniques de construction et th\'eor\`emes
d'existence en g\'eom\'etrie alg\'ebrique IV: Les sch\'emas de Hilbert,}
S\'em. Bourbaki {\bf 13} (1960/61), exp. 221.

\item"{[H1]}" R. Hartshorne, {\it Connectedness of the Hilbert scheme,}
IHES Publ. Math.
{\bf 29} (1966), 261--309.

\item"{[H2]}" R. Hartshorne, {\it Genus of space curves,} preprint 1992 (to
appear in Ann.
Univ. Ferrara, Sez. VII --- Sci. Math.).

\item"{[I]}" A. Iarrobino, {\it Hilbert scheme of points: Overview of last
ten years,}
Proc. Symp. Pure Math. {\bf 46} Part 2 (1987), 297--320.

\item"{[M-D--P]}" M. Martin-Deschamps, D. Perrin, {\it Sur les bornes des
modules de Rao,}
C.R.Acad.Sc. {\bf 317} S\'erie 1 (1993), 1158--1162.

\item"{[Mi]}" Migliore J., {\it On linking double lines,} Trans. A.M.S.
{\bf 294}
(1986), 177-185.

\item"{[P--S]}" R. Piene, M. Schlessinger, {\it On the Hilbert scheme
compactification of the
space of twisted cubics,} Amer. J. Math. {\bf 107} (1985), 761--774.

\item"{[R]}" A. Reeves, {\it Combinatorial structure on the Hilbert
scheme,} Ph. D. Thesis,
Cornell University, 1992.

\item"{[S]}" G. Sacchiero, {\it Fibrati normali de curve razionali dello spazio
proiettivo,} Ann. Univ. Ferrara, Sez. VII --- Sci. Math. {\bf 26} (1980), 33--40.

\item"{[W]}" J. Wahl, {\it Equations defining rational singularities,} Ann.
scient. \'Ec.
Norm. Sup. {\bf 10} (1977), 231--264.
\endroster
\endRefs
\enddocument
\bye